\documentclass[12pt]{article}
\usepackage{amsfonts}
\usepackage{amsmath}
\usepackage{amssymb}

\begin{document}


\renewcommand{\theequation}{\thesection.\arabic{equation}}

\newtheorem{theorem}{Theorem}
\newtheorem{problem}{Problem}
\newtheorem{definition}{Definition}
\newtheorem{lemma}{Lemma}
\renewcommand{\thelemma}{\Alph{lemma}}

\newtheorem{proposition}{Proposition}
\newtheorem{corollary}{Corollary}
\newtheorem{example}{Example}
\newtheorem{conjecture}{Conjecture}
\newtheorem{algorithm}{Algorithm}
\newtheorem{exercise}{Exercise}
\newtheorem{remarkk}{Remark}

\newcommand{\be}{\begin{equation}}
\newcommand{\ee}{\end{equation}}
\newcommand{\bea}{\begin{eqnarray}}
\newcommand{\eea}{\end{eqnarray}}
\newcommand{\beq}[1]{\begin{equation}\label{#1}}
\newcommand{\eeq}{\end{equation}}
\newcommand{\beqn}[1]{\begin{eqnarray}\label{#1}}
\newcommand{\eeqn}{\end{eqnarray}}
\newcommand{\beaa}{\begin{eqnarray*}}
\newcommand{\eeaa}{\end{eqnarray*}}
\newcommand{\req}[1]{(\ref{#1})}

\newcommand{\lip}{\langle}
\newcommand{\rip}{\rangle}
\newcommand{\uu}{\underline}
\newcommand{\oo}{\overline}
\newcommand{\La}{\Lambda}
\newcommand{\la}{\lambda}
\newcommand{\eps}{\varepsilon}
\newcommand{\Om}{\Omega}

\newcommand{\rrr}{{\Bigr)}}
\newcommand{\qqq}{{\Bigl\|}}

\newcommand{\dint}{\displaystyle\int}
\newcommand{\dsum}{\displaystyle\sum}
\newcommand{\dfr}{\displaystyle\frac}
\newcommand{\bige}{\mbox{\Large\it e}}
\newcommand{\integers}{{\mathbb Z}}
\newcommand{\ZZ}{{\Bbb Z}}
\newcommand{\rationals}{{\Bbb Q}}
\newcommand{\reals}{{\Bbb R}}
\newcommand{\realsd}{\reals^d}
\newcommand{\realsn}{\reals^n}
\newcommand{\NN}{{\Bbb N}}
\newcommand{\DD}{{\Bbb D}}
\newcommand{\XX}{{\mathfrak X}}

\newcommand{\degree}{{\scriptscriptstyle \circ }}
\newcommand{\dfn}{\stackrel{\triangle}{=}}
\def\complex{\mathop{\raise .45ex\hbox{${\bf\scriptstyle{|}}$}
     \kern -0.40em {\rm \textstyle{C}}}\nolimits}
\def\hilbert{\mathop{\raise .21ex\hbox{$\bigcirc$}}\kern -1.005em {\rm\textstyle{H}}} 
\newcommand{\RAISE}{{\:\raisebox{.6ex}{$\scriptstyle{>}$}\raisebox{-.3ex}
           {$\scriptstyle{\!\!\!\!\!<}\:$}}} 

\newcommand{\hh}{{\:\raisebox{1.8ex}{$\scriptstyle{\degree}$}\raisebox{.0ex}
           {$\textstyle{\!\!\!\! H}$}}}

\newcommand{\OO}{\won}
\newcommand{\calA}{{\cal A}}
\newcommand{\calB}{{\cal B}}
\newcommand{\calC}{{\cal C}}
\newcommand{\calD}{{\cal D}}
\newcommand{\calE}{{\cal E}}
\newcommand{\calF}{{\cal F}}
\newcommand{\calG}{{\cal G}}
\newcommand{\calH}{{\cal H}}
\newcommand{\calK}{{\cal K}}
\newcommand{\calL}{{\cal L}}
\newcommand{\calM}{{\cal M}}
\newcommand{\calO}{{\cal O}}
\newcommand{\calQ}{{\cal Q}}
\newcommand{\calR}{{\cal R}}
\newcommand{\calX}{{\cal X}}
\newcommand{\calXX}{{\cal X\mbox{\raisebox{.3ex}{$\!\!\!\!\!-$}}}}
\newcommand{\calXXX}{{\cal X\!\!\!\!\!-}}
\newcommand{\gi}{{\raisebox{.0ex}{$\scriptscriptstyle{\cal X}$}
\raisebox{.1ex} {$\scriptstyle{\!\!\!\!-}\:$}}}
\newcommand{\intsim}{\int_0^1\!\!\!\!\!\!\!\!\!\sim}
\newcommand{\intsimt}{\int_0^t\!\!\!\!\!\!\!\!\!\sim}
\newcommand{\pp}{{\partial}}
\newcommand{\al}{{\alpha}}
\newcommand{\sB}{{\cal B}}
\newcommand{\sL}{{\cal L}}
\newcommand{\sF}{{\cal F}}
\newcommand{\sE}{{\cal E}}
\newcommand{\sX}{{\cal X}}
\newcommand{\R}{{\rm I\!R}}
\newcommand{\vp}{\varphi}
\newcommand{\N}{{\rm I\!N}}
\def\ooo{\lip}
\def\ccc{\rip}
\newcommand{\ot}{\hat\otimes}
\newcommand{\rP}{{\Bbb P}}
\newcommand{\bfcdot}{{\mbox{\boldmath$\cdot$}}}

\renewcommand{\varrho}{{\ell}}
\newcommand{\dett}{{\textstyle{\det_2}}}
\newcommand{\sign}{{\mbox{\rm sign}}}
\newcommand{\TE}{{\rm TE}}
\newcommand{\TA}{{\rm TA}}
\newcommand{\E}{{\rm E\,}}
\newcommand{\won}{{\mbox{\bf 1}}}
\newcommand{\Lebn}{{\rm Leb}_n}
\newcommand{\Prob}{{\rm Prob\,}}
\renewcommand{\mod}{{\rm mod\,}}
\newcommand{\sinc}{{\rm sinc\,}}
\newcommand{\ctg}{{\rm ctg\,}}
\newcommand{\loc}{{\rm loc}}
\newcommand{\trace}{{\,\,\rm trace\,\,}}
\newcommand{\Tracenew}{{\,\,\rm trace\,}}
\newcommand{\Dom}{{\rm Dom}}
\newcommand{\ifff}{\mbox{\ if and only if\ }}
\newcommand{\proof}{\noindent {\bf Proof:\ }}
\newcommand{\remark}{\noindent {\bf Remark:\ }}
\newcommand{\remarks}{\noindent {\bf Remarks:\ }}
\newcommand{\note}{\noindent {\bf Note:\ }}
\newcommand{\uuj}{{\uu{j}}}
\newcommand{\uui}{{\uu{i}}}
\newcommand{\on}{{\otimes n}}
\newcommand{\om}{{\otimes m}}

\newcommand{\boldx}{{\bf x}}
\newcommand{\boldX}{{\bf X}}
\newcommand{\boldy}{{\bf y}}
\newcommand{\boldR}{{\bf R}}
\newcommand{\uux}{\uu{x}}
\newcommand{\uuY}{\uu{Y}}

\newcommand{\limn}{\lim_{n \rightarrow \infty}}
\newcommand{\limN}{\lim_{N \rightarrow \infty}}
\newcommand{\limr}{\lim_{r \rightarrow \infty}}
\newcommand{\limd}{\lim_{\delta \rightarrow \infty}}
\newcommand{\limM}{\lim_{M \rightarrow \infty}}
\newcommand{\limsupn}{\limsup_{n \rightarrow \infty}}

\newcommand{\ra}{ \rightarrow }

\newcommand{\ARROW}[1]
  {\begin{array}[t]{c}  \longrightarrow \\[-0.2cm] \textstyle{#1} \end{array} }

\newcommand{\AR}
 {\begin{array}[t]{c}
  \longrightarrow \\[-0.3cm]
  \scriptstyle {n\rightarrow \infty}
  \end{array}}

\newcommand{\pile}[2]
  {\left( \begin{array}{c}  {#1}\\[-0.2cm] {#2} \end{array} \right) }

\newcommand{\floor}[1]{\left\lfloor #1 \right\rfloor}

\newcommand{\mmbox}[1]{\mbox{\scriptsize{#1}}}

\newcommand{\ffrac}[2]
  {\left( \frac{#1}{#2} \right)}

\newcommand{\one}{\frac{1}{n}\:}
\newcommand{\half}{\frac{1}{2}\:}

\def\le{\leq}
\def\ge{\geq}
\def\lt{<}
\def\gt{>}

\def\squarebox#1{\hbox to #1{\hfill\vbox to #1{\vfill}}}
\newcommand{\qed}{\hspace*{\fill}
           \vbox{\hrule\hbox{\vrule\squarebox{.667em}\vrule}\hrule}\bigskip}

\title{Some Ergodic Theorems for  Random Rotations on  Wiener Space}

\author{A. S. \"Ust\"unel and M. Zakai}
\date{}
\maketitle
\begin{abstract}
\noindent
In this paper we study ergodicity and mixing property of some measure
preserving transformations on the Wiener space $(W,H,\mu)$ which are
generated by some random unitary operators defined on the
Cameron-Martin space $H$.
\end{abstract}
\section{Introduction}
Although the Wiener measure is one of the most popular and studied
probability measures, there are surprisingly few results about the
ergodicity of the measure preserving transformations of the Wiener
paths. In fact since the early  work of Maruyama \cite{Ma}  and that  of Wiener
and Akutowicz \cite{W-A}   
it is difficult to find any
work about the subject in the literature. Having studied the general
Wiener-measure  preserving transformations in \cite{U-Z-1} (cf. also
\cite{U-Z-2}), we were led  to look at the ergodicity of these
transformations under the light of the powerful techniques developed
by the Stochastic Calculus of Variations of Paul Malliavin. In fact
even about  the classical Wiener-Ito decomposition we have a better
knowledge now and using this latter technique, one can characterize
rather easily the ergodicity and the mixing property of the
transformations which are the second quantizations of the
deterministic unitary transformations on the Cameron-Martin
space $H$. When  we take a random unitary transformation
$R(w)$ of the Cameron-Martin space with the property that, for any
$h\in H$,  $\nabla Rh$ is a quasi-nilpotent operator on $H$, then $R$
generates also a Wiener-measure preserving transformation, called
rotation. However the ergodicity in this case is much more difficult
to characterize, since the randomness of $R$ induces a very strong
non-linearity and we can not use anymore the Wiener chaos technique as
easily as in the deterministic case. In this paper we study the
ergodicity and the strong mixing property of this kind of random
rotations in some special cases: the first is the case where the
randomness enters as an  input which is independent of the Wiener
``paths''. In the second situation we  replace  this independence hypothesis
by another  interesting hypothesis, namely we assume that the
(Gaussian) divergence of the resolution of the identity associated to
the unitary operator defines a cylindrical martingale (indexed with
the associated spectrum) having the chaotic representation property. 
The next section presents some necessary conditions for ergodicity.
Finally we derive a sufficient condition for strong mixing (which is
also necessary when $R$ is deterministic) and give two generic
examples of strongly mixing transformations. Let us now explain in detail
the plan of the paper.

Let $(W,H, \mu)$ be an abstract Wiener space, i.e. $W$ is a  separable
Banach space, $H$ is a Hilbert space, whose continuous dual is
identified with itself and it is  densely and continuously injected
in $W$. For any $e\in W^*$, 
$$
\int_W e^{i<e,w>}\mu(dw)=E [\exp i<e,w>] = \exp -{|j(e)|_H^2}/{2}\,,
$$
 where $j$ denotes the injection $W^*\hookrightarrow H$. Let $T:W\to W$
 be a measurable, invertible and measure
preserving transformation on $W$, the problem considered in this paper is the
ergodicity of this transformation.  As an example of such a transformation
consider the classical Wiener space: let $B_t, t\in [0,1]$ be a standard
Wiener process taking values in $\realsn$, consider $\gamma(t,w)$ which for
every $t$ in $[0,1]$ takes value in the class of unitary $n\times n$-matrices.
Then, if $\gamma$ is non-random or under suitable measurability assumptions
on $\gamma(t,w)$, the process $Y(t,w)$ defined as
$Y(t,w) = \int_0^t \gamma (\tau,w) dB_\tau$ exists as an Ito integral and,
due to the celebrated theorem of Paul L\'evy, 
$(t,w)\to Y(t,w)$ is also a standard Brownian motion in  $\realsn$.

The class of transformations on Wiener space that will be considered is as
follows.  Let $(e_n,\,n\in \N)$ be a complete, orthonormal basis  of
$H$ and $e_n \in W^*$ for all $n\in \N$  
\footnote{In the sequel, as long as there is no confusion,  we shall
  not distinguish the elements of  $W^*$, from their images in $H$.}.
By the Ito-Nisio theorem (cf. \cite{I-N,U-Z-2})
$$ 
w = \sum_{n=1}^\infty  (\delta e_n(w)) e_n
$$ 
$\mu$-almost surely in the sense that
$$
\left\|w - \sum_1^N (\delta e_n(w)) e_n \right\|_{W} \to 0
$$
$\mu$-almost surely as $N\to\infty$ where $\delta e(w)
=_{W^*}<e,w>_{W}$ is the abstract version of the  Wiener integral of $j(e)$.

Consider first the case where $R$ is a non-random unitary transformation on
$H$, then
$$
w\to T(w) = \sum_{n=1}^\infty \delta(Re_n)(w)\, e_n
$$
is a measure preserving transformation of $W$.
The ergodicity of this class of transformations was characterized in 
\cite{U-Z-4}.  In this paper we consider the problem of
ergodicity for the case where $R=R(w)$ is random. In
\cite{U-Z-1} (cf. also \cite{U-Z-2}) we have already  shown that if $R(w)$ is
almost surely a unitary transformation on $H$, then under additional
(non-trivial) assumptions, the mapping defined by  
\begin{equation}\label{1.1}
w\to T(w) = \sum_{n=1}^\infty \delta \left(R(w) e_n\right)(w) \,\, e_n
\end{equation}
(where $\delta (Re)$ denotes the `divergence' or `Skorohod integral') is a
measure preserving transformation on the Wiener space. Hence the
problem of ergodicity of such transformations is natural.

In the next section we summarize some relevant results from the Malliavin
calculus. In the third section we give necessary and sufficient
conditions for two classes of random rotations in terms of their
resolution of identity. The first class consists of the rotations
whose randomness are independent of the underlying Wiener paths. The
second class maybe described as the set of rotations whose (random)
resolutions  of identity define cylindrical martingales (indexed with
the Cameron-Martin space $H$) with chaotic representation property.

The case of general rotations is considered in Section~4 and necessary
conditions for ergodicity is derived.  A  sufficient condition for
strong  mixing for a general class of rotations is derived in Section~5. 

\section{Preliminaries}

Let $ (W,H,\mu ) $ be an abstract Wiener space, 
a mapping $\varphi$ from $W$ into some separable Hilbert space $\XX$ will be
called a cylindrical function if it is of the form
$ \varphi (w) = f(< v_1 , w>, \cdots, <v_n ,w>) $ where
$f\in C_0^{\infty} (\realsn , \, \XX ) \quad v_i \in W^* $
for $ i= 1\, ,\ldots , \, n $.
For such a $\varphi$, we define $\nabla\varphi$ as
$$
\nabla\varphi (w) = \displaystyle\sum_{i=1}^n\,\partial_i f\Bigl( <v_1 , w>
\, ,\ldots , \, < v_n , w >\Bigr) \, \tilde{v}_i
$$
where $ \tilde{v}_i $ is the image of $ v_i $ under the injection
$ W^* \hookrightarrow H $. It follows that $\nabla$ is a closable operator on
$ L^p (\mu, \XX)$, $ p\ge 1$ and we will denote its closure
with the same notation.
The powers $\nabla^k$ of $\nabla$ are defined by iteration.
For $p\gt 1 $, $ k\ge 1 $, we denote by $ \DD_{p,k} (\XX ) $ the completion of
$\XX$-valued cylindrical functions with respect to the norm:
$$
\|\varphi\|_{\DD_{p,k} (\XX)} \equiv
\|\varphi\|_{p,k} =\displaystyle\sum_{i=0}^k\;
\| \nabla^i\varphi\|_{L^p (\mu,\XX \otimes H^{\otimes i})}
 \quad .
$$
Let us denote by $\delta$ the formal adjoint of $\nabla$ with respect to
the Wiener measure $\mu$ and define $\calL$ as $\delta\circ \nabla $.
The well-known result of P.~A.~Meyer assures that the norm defined above is
equivalent to
$$
|||\varphi|||_{p,k} = \parallel (I+\calL)^{k/2}
\varphi\parallel_{\scriptstyle{L}^p (\mu , \XX)} \quad ,
$$
and  $\calL$ is called  the Ornstein-Uhlenbeck operator or the number
operator.  Note that, due to its self adjointness, its non-integer powers are
well-defined.  Moreover we can also define $ \DD_{p,k} (\XX ) $ for negative
$k$'s using the second norm and we denote by
$\DD (\XX) =
\cap_{p\gt 1}
\cap_{k\in \NN}  \DD_{p,k} (\XX) $ and,
$\DD' (\XX) =\cup_{p\gt 1}
\cup_{k\in\integers}  \DD_{p,k} (\XX) $.
In case $\XX = \reals $ we write simply
$\DD_{p,k} , \, \DD , \, \DD' $ instead of
$\DD_{p,k} (\reals) , \, \DD (\reals) , \, \DD' (\reals) $.
Let us recall that
$$
\nabla \, : \; \DD_{p,k} (\XX)\to \DD_{p, k-1} (\XX\otimes H) 
$$
and
$$ 
\delta \, : \; \DD_{p,k} (\XX\otimes H)\to \DD_{p, k-1} (\XX) 
$$
are continuous linear operators for any $ p\gt 1 , \, k\in\integers $.



\subsection{Rotations}
{\bf Rotation Theorem:}
{\em 
Let $R$ be a strongly measurable random variable on $W$ with values in the
space of bounded linear operators on $H$.  Assume that $R$ is almost
surely  an isometry on $H$(i.e., $|R(w) h|_H= |h|_H$$\mu$-almost
surely,  for any $h\in H$).  Further assume that for some $p>1$ and for all 
$h\in H$, $Rh \in \DD_{p,2}(H)$ and $\nabla Rh \in \DD_{p,1} (H\otimes H)$
is a quasi-nilpotent operator on $H$\footnote{This means that 
$\lim_{n\to\infty} \|(\nabla Rh)^n\|_{L(H,H)}^{1/n}= 0$ almost surely  or,
equivalently, $\trace (\nabla Rh)^n = 0$ almost surely, for all
$n \ge 2$.}. If moreover, either
\begin{itemize}
\item[a)] for any $h\in H$, 
$$
(I+i\nabla Rh)^{-1} \cdot Rh\in L^q(\mu,H),\, q > 1
$$  
(here $q$ may depend on $h\in H$) or,
\item[b)]~$Rh \in \DD(H)$ for any $h\in H$, 
\end{itemize}
then
\begin{equation}
\label{2.1new}
E\Bigl[\exp i \delta (Rh)\Bigr] = e^{ - \half |h|^2_H}\,.
\end{equation}
Besides, for any complete, orthonormal basis $(e_i,i\geq 1)$ of $H$,
the sum
$$
T(w)=\sum_{i=1}^\infty \delta(Re_i)(w)e_i
$$
converges almost surely in the strong topology of $W$, the result is
almost surely independent of any particular choice of $(e_i,i\geq 1)$,
consequently $T$ defines a measure preserving transformation of $W$
which is called the rotation associated to $R$.
}

\remark In fact it suffices to assume (\cite{U-Z-2}) the above
hypothesis for any 
$h\in H_1$, where $H_1$ is any arbitrary dense vector subspace of
$H$.\\

\noindent
>From this theorem it follows that
$H$, 
$$
(\delta(Re_i), i\geq 1)
$$ 
are independent, identically distributed
(i.i.d.) $N(0,1)$-random variables and the  equation \req{1.1} defines 
a measure preserving transformation of $W$ thanks to the Ito-Nisio
theorem (cf.\cite{I-N,U-Z-2}).
The random isometry  $R$ satisfying the conditions for this theorem
under (b) will be said to
{\bf  satisfy the rotation conditions}.
Let us remark that, to  an  operator $R$ with the  above properties, for any 
fixed $k\in H$, it corresponds another one, satisfying the same properties,  
defined as $w\to R(w+tk)$, 
$t\in [0,1]$, that we shall denote by $R_{t,k}$. With this notion we
define a new operator as  
$$
X_k^RF(w)=\frac{d}{dt}F(T_{t,k}(w))|_{t=0}\,,
$$
where $T_{t,k}w$ is defined as 
$$
T_{t,k}w=\sum_{i=1}^\infty\delta(R_{t,k}e_i)(w)e_i\,.
$$
$X^R$ is closable (\cite{U-Z-2})  and we have 
\begin{equation}
\label{total-der}
\nabla_k(F\circ T)=(R(\nabla F\circ T),k)_H+X_k^RF
\end{equation}
for any cylindrical $F$. This operator plays an important role in the
analysis  of random rotations:
\begin{lemma}
\label{compo-lemma}
Let $u:W\to H$ be any cylindrical map, then one has 
$$
\delta u\circ T=\delta(R(u\circ T))+ \trace(R X^Ru)\,.
$$
\end{lemma}
\proof
Let $(e_i,i\in \NN)$ be a complete, orthonormal basis of $H$. 
We have, using the relation (\ref{total-der}) and denoting $(u,e_i)_H$ by 
$u_i$, 
\begin{eqnarray*}
\delta u\circ T&=&\sum_{i=1}^\infty\left\{u_i\circ T\,\delta(Re_i)-
((\nabla_{e_i}u_i)\circ T\right\}\\
&=&\sum_{i=1}^\infty\left\{u_i\circ T\,\delta(Re_i)-
(R\nabla u_i\circ T,Re_i)_H\right\}\\
&=&\sum_{i=1}^\infty\left\{u_i\circ T\,\delta(Re_i)-
(\nabla(u_i\circ T)-X^Ru_i,Re_i)_H\right\}\\
&=&\delta(R(u\circ T))+\sum_{i=1}^\infty(X^Ru_i,Re_i)_H\,.
\end{eqnarray*}
\qed

\remark
Since $\delta u\circ T$ and $\delta(R(u\circ T))$ are independent of the 
choice of $(e_i,i\in \NN)$, so does $\trace(R X^R u)$.

\subsection{Traces}
 Let $H$ be a separable Hilbert space and let
$\varphi =( \varphi_i,\,i\geq 1)$ be a {\em fixed\/}
complete orthonormal basis  of  $H$.
We will use $\varphi^{(n)} =(\varphi_\uui^{(n)})$ to denote the
complete orthonormal basis on $H^\on$ induced by $\varphi$,
i.e.\
$\uui = (i_i, \ldots, i_n)\in \NN^n$,
$\varphi_\uui^{(n)} = \varphi_{i_1} \otimes \varphi_{i_2} 
\otimes \cdots \otimes \varphi_{i_n} $
and the sequence $\varphi_\uui^{(n)}$'s  are 
arranged in lexicographical order.  For
a bounded operator $A$ on $H^\on$ we define the $\varphi$-trace as:
\begin{equation}\label{2.1}
\trace^\varphi A = \sum_\uui ( A \varphi_\uui^{(n)},\varphi_\uui^{(n)})_{H^\on}
\end{equation}
where the summation is in lexicographical order and provided the series
converges.  From now on we will delete the $\varphi$ and denote the series
defined by \req{2.1} as $\trace\,A$.

Let $ u\in\DD_{2,1}(H) $, 
the $\varphi$-Ogawa integral of $u$ is defined as
$$ 
\delta^{\varphi} \circ u =\dsum_i \, (u , \, \varphi_i) \, \delta \varphi_i
$$
provided that the series converges in $ L^2 $.
Then \cite{Z}, $ \delta^{\varphi} \circ u $ exists iff
$ \Tracenew^{\varphi}\, \nabla u $ 
exists and then
$ \delta u=\delta^{\varphi} \circ u -\Tracenew^{\varphi}\, \nabla u $.


The following two lemmas will be needed later.

\begin{lemma}
Let $R$ satisfy the rotation conditions, and for some fixed
complete orthonormal basis $\varphi=(\varphi_i,\,i\geq 1)$,
$\mu$-almost surely 
$$
\trace^\varphi \nabla R(w) h=0, 
$$
for any $h\in H$. 
Then
$$
\sum_i (\delta \varphi_i) \,\, R^* \varphi_i 
= \sum_i \delta (R \varphi_i) \,\, \varphi_i\,.
$$
\end{lemma}
\remark Note that the right hand side  is independent of $\varphi$,
but equality holds only if $\trace^\varphi \nabla Rh=0$.\\
\proof Let $h$ be an element of $H$, then 
\beaa
\sum \delta \varphi_i (R^* \varphi_i, h) & = &
\sum_i \delta \varphi_i (\varphi_i, Rh)_H \\
& = & \delta^\varphi \circ (Rh) \\
& = & \delta (Rh)  \\
& = & \sum_i (\delta R \varphi_i) \cdot (\varphi_i, h)_H\,.
\eeaa
\qed

\begin{lemma}
Let $R$ and $\varphi$ be as in the lemma above, and
let $v \in \DD(H)$, 
$\trace^\varphi \nabla v=0$ and 
$\trace^\varphi \nabla (R(w) v(Tw)) = 0$. 
Then
$$
( \delta v )\circ T = \delta 
\Bigl( R (v\circ T) \Bigr)\,.
$$
If only $\trace^\varphi \nabla v =0$ then
\begin{equation}
\label{new2.5}
(\delta v) \circ T = \delta \circ \Bigl(R( v\circ T) \Bigr)
\end{equation}
\end{lemma}

\proof
Note that 
\begin{equation}
\label{*}
\delta h \circ T = \sum \delta (R\varphi_i)(\varphi_i, h)_H
= \delta Rh\,.
\end{equation}

Now, by \req{*} and Lemma~A,
\beaa
\delta v \circ T &=&
(\delta^\varphi \circ v + \trace^\varphi \nabla v) \circ T \\
& = & (\delta^\varphi \circ v) \circ T \\
& = & \sum\left((v, \varphi_i)_H \delta \varphi_i\right) \circ T\\
& = & \sum\left(v\circ T, \varphi_i \right)_H \delta R \varphi_i \\
& = & \sum \left(v\circ T,  R^* \varphi_i\right)_H \delta \varphi_i \\
& = & \delta \circ \left(R( v\circ T) \right)\hspace{3cm} 
\mbox{(this proves \req{new2.5})} \\
& = & \delta \left(R( v\circ T) \right)
\eeaa
\qed

\section{Chaos representation and ergodic  rotations}
\setcounter{equation}{0}

Let $(W , H, \mu)$ be an  abstract Wiener space.
Let $(p_\theta,\,\theta\in [0, 2\pi])$ be a right continuous resolution of
 identity on $H$ and let $\calR$ denote the class of non random 
unitary operators
 on $H$ which are represented
by it:
$$
\calR = \left\{ R: \, R = \int_0^{2\pi} e^{ i 
\varphi (\theta)} d p_\theta \right\} 
$$
where $\varphi (\cdot)$ is real valued,  right continuous on 
$[0, 2\pi]$.  Note that the elements of $\calR$ commute.
Further assume that $\int_A d p_\theta\not= 0$ for all $\theta$ sets $A$ of
 positive Lebesgue measure.

Let $(M,\calM, P)$ be a probability space,
independent of $W$.
Let $(R_i(m), i\in \NN)$ be an i.i.d. sequence taking 
values in $\calR$ and $(\psi_i(\theta,m), i\geq 1) $ 
are $\calM$-measurable i.i.d.\ continuous
functions on $[0, 2\pi]$.  Consider the product space
$(W\times M, \calB(W)\otimes \calM, \mu \times P)$.
Set
$$
R_n(m) = \int_0^{2\pi} e^{ i \psi_n (\theta,m)} dp_\theta\,.
$$
Let us define $T^i(w,m)$ as
$$
T^i(w,m)=\sum_{i=1}^\infty \delta\left(R_1(m)R_2(m)\ldots
  R_n(m)e_i\right)e_i\,.
$$
By ergodicity we mean that for any square integrable $F$, 
$$
\lim_{n\to \infty}\one \sum_{i=1}^n F(T^i(w,m))= E[F]
$$
$\mu\times P$-almost surely.
Recall that from  Lemma~A
$$
T^1(w,m) = \sum_i \delta e_i \cdot R_1^* (m) e_i 
$$
and
$$
T^n(w,m) = \sum_i (\delta e_i)\, R_n^*
(m)\, R_{n-1}^* (m) \cdots R_1^* (m) e_i\,.
$$

\begin{theorem}
\label{chaos-thm}
Under the above assumptions on $R$ and $T$,  the following is
necessary and sufficient for the ergodicity of $T=(T^i, i\in\NN)$ :
\begin{enumerate}
\item $\theta\to (p_\theta h, h)_H$ is continuous on $[0, 2\pi]$ for all
$h\in H$.
\item For all $\eta$ in $[0, 2\pi]$, the inequality
\begin{equation}
\Bigl|E \left[\exp i \Bigl( \psi_1 (\theta,m) - \eta \Bigr)\right] \Bigr|
< 1
\end{equation}
holds for almost all  $\theta$ in $[0, 2 \pi]$ with respect to the
Lebesgue measure.
\end{enumerate}
\end{theorem}

\proof
Starting with necessity, assume that $\theta\to p_\theta$ is discontinuous at
$\theta=\theta_0$. If $h$ is in the invariant subspace defined by the
projection $p_{\theta_0} - p_{\theta_0-}$,
then
$$
R_n(m)  h = e^{ i \psi_{n} (\theta_0,m)}  h
$$
for any $n\geq 1$, hence, for a.a.\ $m$
$$
| \delta R_n \, h| = |\delta h |
$$
and $T$ is not ergodic.  Similarly, assume that
$$
E\left[e^{ i \psi_1 (\theta, m)}\right] = e^{i \eta}
$$
holds for almost all $\theta\in A$, where $A$ is a measurable subset
of $[0,2\pi]$ of positive Lebesgue measure. Then
$\psi_1(\theta,m)=\eta$ almost surely on $A\times M$. 
Set $H_A$ to be the invariant subspace of the projection
$\int_A d p_\theta$. For any  $h\in H_A$ and $n\geq 1$, we have
$$
\delta h = e^{-i \eta} \delta R_nh
$$
and, again, this result implies that  $T$ is non ergodic.

Turning now to sufficiency, assume that $T$ is not ergodic, then for some 
square integrable $F$,
\begin{equation}
\label{inva-eq}
F(w,m) = F(T^1(w, m))\qquad \mbox{a.s.\ } \mu \times P\,,
\end{equation}
fixing $m$ and developing $F$ in a multiple Wiener-Ito series:
$$
F=\sum_{n=0}^\infty I_n(K_n)\,,
$$
the relation (\ref{inva-eq}) yields 
$$
I_n (K_n^{(m)}) = I_n (R_1^{\otimes n} (m) K_n^{(m)}) 
$$
$ \mu\times P$-almost surely. Hence, $P$-almost surely
\begin{eqnarray}
\label{ostar}
\lefteqn{
0= \int_{[0, 2\pi]^n}
\left|1- \prod_{j=1}^n \exp i \psi_1 (\theta_j, m) \right|^2 
d\Bigl( ( p_{\theta_1} \otimes \ldots \otimes p_{\theta_n})
K_n^{(m)}, K_n^{(m)} \Bigr)_{H^{\odot n}}}
\nonumber \\
&=&
\int_{[0,2\pi]^n} \left|
1-\prod_{j=1}^n \exp i \psi_1 (\theta_j, m) \right|^2
\rho (d \theta_1, \cdots , d\theta_n, m)
\end{eqnarray}
where $\rho$ is a continuous (atomless)  positive measure
(cf. \cite{U-Z-4}). Hence
$$
0= \int_{[0, 2\pi]^n}
\left|1-\cos \sum_{j=1}^n \psi_1 (\theta_j, m)\right|
\rho (d\theta_1, \ldots , d\theta_n,m)\,.
$$
Consequently 
\begin{equation}
\label{conc-con}
\sum_{j=1}^n \psi_1 (\theta_j, m) = 0\; \mod 2 \pi
\end{equation}
$P$-almost surely and the support of  $\sum_{j=1}^n \psi_1 (\theta_j,
m)$ lies in a sub-manifold of
$[0, 2\pi]^n$ whose dimension is at most $n-1$.  Since the measure
$\rho$ is continuous,  it vanishes on the  lower 
dimensional manifolds. Consequently (\ref{conc-con}) is impossible and $T$
is ergodic. 
\qed

The following result is almost a corollary of Theorem \ref{chaos-thm}:
\begin{theorem}
\label{chaos-2-thm}
Let  $R$ be  a weakly measurable random variable with values in the
set of unitary operators on $H$ satisfying the rotation
condition. Assume that it has a  representation as  
$$
R=\int_{[0,2\pi]}e^{i\theta}dp_\theta (w)\,,
$$
where $(p_\theta (w),\theta\in [0,2\pi])$ is a weakly measurable
resolution of identity on $H$. Assume furthermore that, for any $h\in
H$,  
$$
\theta \to \delta p_\theta h=m_\theta(h)
$$
is a martingale with respect to the filtration
$(\calD_\theta,\,\theta\in [0,2\pi])$, whose predictable increasing
process, denoted by $(a_\theta(h,h),\,\theta\in [0,2\pi])$ is
deterministic, where $\calD_\theta$ denotes the 
right continuous filtration generated by $\{\delta p_\tau h,\,h\in
H,\,\tau\leq \theta\}$. Then the transformation $T:W\to W$ is ergodic
if and only if the vector measure defined by $\theta\to
d_\theta E[(p_\theta h,h)_H]$ has no atom.
\end{theorem}
\proof 
Let us note that, since $\delta h=m_{2\pi}(h)$, the cylindrical
martingale $(m_\theta,\theta\in [0,2\pi])$ has the chaotic
representation property. Besides, we have 
\begin{eqnarray*}
E[m_\theta(h)^2]&=&E[E[\delta h|\calD_\theta]^2]\\
               &=&E[\delta h\,\delta p_\theta h]\\
              &=&E[(h,p_\theta h)_H]\,,
\end{eqnarray*}
hence $a_\theta(h,h)=E[(p_\theta h,h)_H]$ for any $h\in H$ and
$\theta\in [0,2\pi]$. 
The chaotic representation property means that any square integrable random
variable $F$ can be represented as 
$$
F=E[F]+\sum_{n=1}^\infty
\int_{[0,2\pi]^n}(f_n(t_1,\ldots,t_n),dm_{t_1}\otimes\cdots\otimes
dm_{t_n})_{H^{\otimes n}}\,,
$$
where $f_n:[0,2\pi]^n\to H^{\otimes n}$ is measurable, symmetric with
respect to $(t_1,\ldots, t_n)$ and 
$$
E|F-E[F]|^2=\sum_{n=1}^\infty n! \int_{[0,2\pi]^n}d((a_{t_1}\otimes\cdots
\otimes a_{t_n})f_n,f_n)_{H^{\otimes n}}\,.
$$ 
Note that 
$$
F\circ T=E[F]+\sum_{n=1}^\infty
\int_{[0,2\pi]^n}e^{i\sum_{k=1}^nt_k}
(f_n(t_1,\cdots,t_n),dm_{t_1}\otimes\ldots\otimes dm_{t_n})_{H^{\otimes n}}\,,
$$
hence the rest of the proof goes exactly as the proof of Theorem
\ref{chaos-thm}. 
\qed

\section{A necessary condition for ergodicity of
\newline non-independent rotations}
\setcounter{equation}{0}

Let $R:W\to O(H)$ be a random unitary operator satisfying the rotation
condition. Suppose that it has a representation given as 
\begin{equation}
\label{4.1}
R(w) = \int_0^{2\pi} e^{i \psi (\theta, w)} dp_\theta\,,
\end{equation}
where  the random function $\psi$ takes values in the class of
continuous Lebesgue measurable functions from $[0,2\pi]$ to $[0,2\pi]$.   
\begin{proposition}
\label{mix-prop}
Assume that  $\psi(\theta, w) \in
\DD_{2,1}(L^2([0,2\pi],d\theta))$ and that  $\nabla \psi (\theta, w)$ is
orthogonal to the subspace 
induced by $p_\theta - p_{\theta_-} $ for every $\theta$.  Then the
following  conditions
\begin{itemize}
\item[(a)]
$p_\theta$ is continuous on $[0,2\pi]$
\item[(b)] If $A$ is nonrandom Lebesgue
measurable subset of $[0,2\pi]$ such that for some $\eta$ and for
a.a. $\theta \in A$, 
\quad
$ \psi (\theta, w) = \eta$  almost surely \ 
then the Lebesgue measure of $A$ is zero.
\end{itemize}
are necessary  for the  ergodicity of $T$ which is generated by $R$.
\end{proposition}

\remarks
\begin{enumerate}
\item Equation~\req{4.1} implies that $R$ and $R\circ T$ commute.  
\item The requirement that $\nabla \psi (\theta, w)$ is orthogonal to
$p_\theta-p_{\theta-}$ is satisfied if $\psi (\theta,w)$ is predictable
with respect  to
the $\sigma$-field generated  by 
$\{\delta(p_\theta h), h \in H\}$. 
\item Condition (b) is a necessary condition under
\req{4.1} even if the orthogonality condition for $p_\theta-p_{\theta_-}$ 
is not satisfied.
\end{enumerate}
\proof
Assume that $p_\theta$ is
discontinuous at $\theta=\theta_0$.  Then there exists $h\in H$ such that
$(p_{\theta_0} - p_{\theta_0-}) h = h$ and then
\begin{align*}
\delta (Rh) & = \delta (e^{i\psi(\theta_0,w)} h) \\
& = e^{i\psi (\theta_0, w)} \delta (h)
\end{align*}
by the orthogonality assumption for $\nabla \psi$.  Hence
$$
|\delta Rh| = | \delta h|
$$
and $|\delta h|$ is a nontrivial eigenfunction with eigenvalue 1,
therefore $T$ can not be ergodic.  Similarly, assume that for some
$\eta$ and a set $A$ of positive Lebesgue measure
$$
\exp i \psi (\theta, w) = e^{i \eta}
$$
for a.a.\ $\theta \in A$, then for $\pi_A = \int_A d_\theta p_\theta$ and $h$
invariant with respect to  $\pi_A$
$$
\delta Rh = e^{i\eta} \delta h
$$
and again $|\delta h|$ is invariant hence $T$ is  not ergodic.
\qed
\section{A condition for mixing}
\setcounter{equation}{0}
In this section we give a  sufficient condition for the
strong mixing property of some random rotations.
\begin{theorem}
\label{mix-thm}
Let $R$ satisfy the rotation condition, define inductively the
sequence of operators $(Q_n,\,n\geq 1)$ as  $Q_1(w) = R(w)$ and 
$$
Q_n (w) = R(w) \cdot R(Tw) \cdots R(T^{n-1} w)
$$
for $n\geq 2$.
Assume that for  all $n\in \NN$ and all $k,\,h\in H$
\begin{equation}
\label{star}
\delta h\circ T^n=\delta(Q_nh)
\end{equation}
almost surely and the random variable
\begin{equation}
\label{star-1}
w\to \delta \left(Q_n(\cdot+k)h\right)(w)
\end{equation}
has  a  Gaussian distribution  with variance $|h|_H^2$. 
Then $T$ is strongly  mixing  if, for any $h,\,k\in H$, 
\begin{equation}
\label{5.2}
\lim_{n\to \infty}(k,Q_n h)_H= 0 
\end{equation}
in probability.
\end{theorem}
\begin{remarkk}{\rm
 Before the proof of the theorem, let us give  some typical examples of
 situations in which the conditions (\ref{star}) and (\ref{star-1}) hold:
\begin{enumerate}
\item
By Lemma~B (Section~2), if for all $ n\in \NN$ and $ h\in H$ 
\begin{equation}
\label{starstar}
\trace^\varphi \nabla Q_n h =0
\end{equation}
then condition (\ref{star}) holds.
\item If, for any $h\in H$, $\nabla Q_n h$ is quasi-nilpotent or if
  $Q_nh$ is adapted to 
  the standard Wiener filtration, then the condition (\ref{star-1}) holds.
\item Assume that $W=C([0,1],\reals^d)$ and that  $\sigma:[0,1]\times W\to
  O(\reals^d)$ (orthogonal transformations of $\reals^d$) is an
  optional  process. Define $R$ as 
$$
R(w)h(t)=\int_0^t\sigma(s,w)\dot{h}(s)ds,\,\,h\in H\,.
$$
Then the transformation $T$ defined as
$$
T(w)=\sum_{i=1}^\infty \delta(Re_i)\,e_i
$$
satisfies  the hypothesis (\ref{star}) and (\ref{star-1}) (cf.\ also
the example at the end of this section).
\end{enumerate}
\noindent
In fact to see the last claim assume that $\dot{u}$ is a $dt\times
d\mu$-square integrable, smooth  optional step process and let $u$ be the
$H$-valued random variable whose Lebesgue density is $\dot{u}$. Then we have 
\begin{eqnarray*}
\delta u\circ T&=&\sum_i(\dot{u}_{s_i},W_{s_{i+1}}-W_{s_i})\circ T\\
&=&\sum_i(\dot{u}_{s_i}\circ T,W_{s_{i+1}}\circ T-W_{s_i}\circ T)\\
&=&\sum_i\left(\dot{u}_{s_i}\circ T,\delta( RU_{[s_i,s_{i+1}]})\right)\,,
\end{eqnarray*}
where $U_{[s_i,s_{i+1}]}$ denotes the image in $H$ of the indicator
function of the interval $[s_i,s_{i+1}]$ under the usual injection of
$L^2([0,1])$ into $H$, i.e. $f(s)\to \int_0^\cdot f(s)ds$ and
$(W_t,t\in [0,1])$ is the $d$-dimensional Wiener process.  We also
have from Lemma \ref{compo-lemma}
\begin{eqnarray}
\label{calcul}
\lefteqn{\left(RU_{[s_i,s_{i+1}]},\nabla(\dot{u}_{s_i}\circ
    T)\right)}\nonumber\\ 
&=&(RU_{[s_i,s_{i+1}]},R\nabla\dot{u}_{s_i}\circ T)+
(RU_{[s_i,s_{i+1}]},X^R\dot{u}_{s_i})\nonumber\\
&=&(U_{[s_i,s_{i+1}]},\nabla\dot{u}_{s_i}\circ T)+
(RU_{[s_i,s_{i+1}]},X^R\dot{u}_{s_i})\\
&=&0\,,\nonumber 
\end{eqnarray}
where the first term at (\ref{calcul}) is zero because $\dot{u}_{s_i}$
is $\calF_{s_i}$-measurable, hence its derivative has its support in
the interval $[0,s_i]$. For the second term, it suffices to take
$\dot{u}_{s_i}$ of the form $f(\delta l)$, where $f$ is a smooth
function, $l\in H$ such that the support of $\dot{l}$ is in
$[0,s_i]$. Then we have 
\begin{eqnarray*}
(RU_{[s_i,s_{i+1}]},X^R\dot{u}_{s_i})&=&f'(\delta Rl)(\delta(\nabla
Rl),RU_{[s_i,s_{i+1}]})\\ 
&=&f'(\delta Rl)\delta(\nabla_{RU_{[s_i,s_{i+1}]}}Rl)
\end{eqnarray*}
and it is immediate to see that $\nabla_{RU_{[s_i,s_{i+1}]}}Rl=0$
because of the special form of $R$.
Hence, we see that $(\delta u)\circ T=\delta (R(u\circ T))$, then the
general case follows by a limiting   argument.}
\end{remarkk}
{\bf{Proof of the theorem:}}
We will show that $\lim_{n\to \infty}E[F\,\, F\circ T^n]= 0$ for all
square integrable $F$ 
such that $E[F]=0$ and this implies mixing.  Since the span of the Wick
exponentials is dense in $L^2(\mu)$, it suffices to show that
\begin{equation}
\label{5.4}
E\left[ \rho (\delta k) \rho \Bigl( \delta (Q_n h)\Bigr)\right] \AR 1
\end{equation}
for all $ h,k \in H$, where
$ \rho(\delta k) =\exp (\delta k - \half |k|^2_H)$
and
\begin{eqnarray*}
\rho(\delta (Q_n h))&=& \exp\left\{\delta (Q_n h) - \half |h|_H^2\right\}\\
&=& \rho (\delta h) \circ T^n\,.
\end{eqnarray*}   
Again by a density argument, it suffices to show that
$$
E\left[(\delta k)^l \rho \left( \delta (Q_n h) \right)\right]
\AR E[(\delta k)^l]\,,
$$
for any $l\in \NN$.
 By Theorem 3.5.4 and  Corollary 3.6.1 of \cite{U-Z-2}
\beaa
E[(\delta k)^l] &=& E \left[(\delta k)^l (w + Q_n h) \rho (-\delta Q_n h
) \right] \\
& =& E \left[\left(\delta k + (k, Q_n h)_H\right)^l \cdot \rho
  (-\delta Q_n h)\right] \,.
\eeaa
Since $(k, Q_n h)_H$ is bounded and converges to zero in probability
\beaa
\limn E\left[ (\delta k)^l \rho (-\delta(Q_n h))\right] &=&
\limn E\left[ \left(\delta k + (k, Q_n h)_H\right)^l \, \rho (-\delta Q_n h)
\right]\\
&=& E[(\delta k)^l]\,.
\eeaa
\qed

\remark
Note that the condition (\ref{5.2}) is also necessary when $R$ is a
deterministic operator. More generally, if $T$ is strongly mixing in
the frame of a classical Wiener space, the Ito representation theorem
implies that 
$$
\lim_{n\to \infty}(Q_nh,k)_H=0
$$
in the weak $L^p$-topology for any $p\geq 1$.

An example for a rotation satisfying condition \req{star} (via
\req{starstar}) is the following:\\
Assume that
\begin{itemize}
\item[(a)]
$R(w) = \int_0^{2\pi} e^{ i \psi (\theta, w)} d p_\theta$
and
\item[(b)]
$\psi(\theta, w)$ is adapted to $\calF_\theta = \sigma
\{\delta(p_\theta h), h\in H\}$,
then $\psi(\theta, Tw)$ is also $\calF_\theta$ adapted.  
\end{itemize}
Now,
\begin{equation}
\label{newstar}
Q_n(w) = \int_0^{2\pi} \exp \sum_{j=1}^n i\psi (\theta, T^j w) d p_\theta
\,.
\end{equation}
Then
$$
\nabla Q_n (w) h = 
\sum_{k=1}^n \int_0^{2\pi} i \exp i \sum_{j=1}^n \psi (\theta, T^jw)
\nabla \psi (\theta, T^k w) d p_\theta\,.
$$
Under suitable smoothness conditions and since $\psi$ is adapted it holds that
$$
\nabla \psi (\theta, T^k w)\; \bot\; (p_{\theta_2} - p_{\theta_1}) h,
\,, \mbox{\ if \ }
\theta_2 > \theta_1  \ge \theta\,,
$$
hence $\trace^\varphi \nabla Q_n(w) h = 0$ and
\req{star} also holds.

This result  can be generalized as the following theorem, the proof of
which goes exactly along the same lines as the proof of Theorem
\ref{mix-thm}, hence it will be  omitted:
\begin{theorem}
\label{mix2-thm}
Assume that $(Q_n,n\geq 1)$ is a sequence of random isometries of $H$
such that $Q_nh$ is in the domain of the divergence operator and
$\delta(Q_nh)$ is an $N_1(0,|h|_H^2)$-Gaussian random variable for any
$h\in H$. Assume moreover that the shift defined as $w\to w+Q_n(w)h$
satisfies the Girsanov identity, in the sense that 
$$
E\left[F(w+Q_n(w)h)\exp\left\{-\delta(Q_nh)-\frac{1}{2}|h|_H^2\right\}\right]=
E[F]
$$
for any $F\in C_b(W)$.  Denote by $T_n$ the measure preserving
transformation of $W$, defined by $Q_n$, i.e. $\delta h\circ
T_n=\delta(Q_nh)$, $h\in H$. Then a sufficient condition for the
strong mixing property of $(T_n,n\geq 1)$ is that 
$$
\lim_{n\to \infty}(Q_nh,k)_H=0
$$
in probability, for any $h,\,k\in H$.
\end{theorem}
Here is an application of Theorem \ref{mix-thm}:
\begin{example}
Let $W=C_0([0,1],\R^d)$, then the Cameron-Martin space is the space of
the $\R^d$-valued, absolutely continuous functions on $[0,1]$, with
the square integrable derivatives. Assume that $R$ is given by 
$$
Rh(t)=\int_0^t R_t h'(t)dt\,,
$$
where $R_t$ is an $\R^d\otimes \R^d$-valued, adapted process such
that, for any $x\in \R^d$, $|R_t^\star R_tx|=|x|$ almost
surely. Define $T$ as to be $\delta h\circ T=\delta(Rh),\,h\in
H$. Assume that $R_t\otimes R_s$ is independent of $(R_t\otimes
R_s)\circ T\ldots (R_t\otimes R_s)\circ T^{n-1}$ for any $n\geq 2$,
$s<t\in [0,1]$-$ds\times dt$ almost surely. Assume moreover that the two point
function $A_{s,t}=E[R_s\otimes R_t]$ satisfies the following:
$$
\lim_{n\to \infty}(A_{s,t}^nx,y)_{\R^{2d}}=0
$$
almost surely for any $x,y\in \R^{2d}$, $s<t\in [0,1]$. Then $T$ is
strongly mixing. 
\end{example}
Let us give another  example:
\begin{example}
Assume that $W=C_0([0,1],\R)$, with the corresponding Cameron-Martin
space. Assume also that $((b^i_t,t\in [0,1]),\,i\geq 1)$ is a sequence
of one-dimensional Wiener processes, independent of $W$. Define
$(T_n,n\geq 1)$ inductively  as 
\begin{eqnarray*}
T_1w(t)&=&w_1(t)=\int_0^t\sign(b^1_s)dw_s\,,\\
T_{n+1}w(t)&=&w_{n+1}=\int_0^t\sign(b^n_s)dw_n(s)
\end{eqnarray*}
and regard $T_n$ as a function of the Wiener path although it depends
also on $b^1,\ldots,b^n$. Then it is a measure preserving
transformation of $W$. We have 
\begin{eqnarray*}
\delta h\circ T_n&=&\delta Q_nh\\
        &=&\int_0^1\sign(b^n_s)\ldots\sign(b^1_s)h'(s)dw_s\,.
\end{eqnarray*}
Since $ds\times dt$-almost surely,  $|E[\sign(b^i_s)\sign(b^i_t)]|<1$,
we have, for any $h,k\in H$,  $\lim_{n\to \infty}(Q_nh,k)_H= 0$ in
$L^2$, hence the sequence $(T_n,n\geq 1)$ is strongly mixing.

\end{example}


{\footnotesize}
\begin{tabular}{ll}
A.S. \"Ust\"unel, & M. Zakai,\\
ENST,   D\'ept. R\'eseaux, \hspace{1.5cm }
& Department of Electrical Engineering,\\
46 Rue Barrault, & Technion---Israel Institute of Technology,  \\
75013 Paris&Haifa 32000,  \\
France&Israel\\
ustunel@enst.fr & zakai
\end{tabular}
\end{document}